\documentclass[12pt]{amsart}

\usepackage[english]{babel}
\usepackage{amsmath,amsfonts,amssymb,amscd}
\usepackage[shortlabels]{enumitem}
\usepackage{geometry}
\usepackage{microtype}
\usepackage[colorlinks=true,citecolor=blue,linkcolor=blue,urlcolor=blue]{hyperref}
\usepackage[capitalise,noabbrev]{cleveref}

\newtheorem{theorem}{Theorem}[section]
\newtheorem{corollary}[theorem]{Corollary}
\newtheorem{lemma}[theorem]{Lemma}
\newtheorem{proposition}[theorem]{Proposition}
\newtheorem{conjecture}{Conjecture}[section]
\newtheorem{question}{Question}[section]
\theoremstyle{definition}
\newtheorem{definition}{Definition}[section]
\newtheorem{example}{Example}[section]
\theoremstyle{remark}
\newtheorem{remark}{Remark}[section]

\newcommand{\Z}{\mathcal Z}
\newcommand{\Wh}{\operatorname{Wh}}
\newcommand{\Fr}{\operatorname{Fr}}
\newcommand{\Int}{\operatorname{Int}}
\newcommand{\id}{\operatorname{id}}
\newcommand{\propi}{\operatorname{pro}\text{-}\pi_1}

\begin{document}

\title[On $\Z$-compactifiability of manifolds]
{On $\Z$-compactifiability of manifolds}

\author{Shijie Gu}
\address{Department of Mathematics\\
Northeastern University, Shenyang, Liaoning, China, 110004}
\email{shijiegutop@gmail.com}
\thanks{This research was supported in part by NSFC grant 12201102.}
\date{Revised August 24, 2026}

\begin{abstract}
In 1976, Chapman and Siebenmann established necessary and sufficient
conditions for $\Z$-compactifying Hilbert cube manifolds.  Although the
corresponding conditions are known to be necessary for a manifold $M^n$ to
admit a $\Z$-compactification, it remains open whether they are sufficient.
Guilbault and the author proved that they are sufficient for
$M^n\times I$, when $n\geq5$.  We further explore this question by giving
additional hypotheses under which the interval factor can be removed.  A
retraction defined near the central added set is sufficient; a
product-compatible splitting also gives a collar; and an
upper-semicontinuous cell-like decomposition gives a splitting-free
criterion.  We also answer affirmatively a question of Guilbault--Tinsley by showing that: for every
$n\ge6$ there is a connected one-ended open PL $n$-manifold
which is $\Z$-compactifiable but not pseudo-collarable. Finally, we discuss the additional control needed
for applications to universal covers of closed aspherical manifolds.
\end{abstract}

\maketitle

\section{Introduction}

This paper sets forth technical variations on its predecessors
\cite{GG20,Gu20} in order to study the $\Z$-compactifiability of manifolds
and further clarify its relationship with pseudo-collarability.  These
notions arise from the following kind of ``nice'' compactification.  A
topological (TOP) $n$-manifold $M^n$, possibly with boundary, is
\emph{completable} if there is a compact topological manifold
$\widehat M^n$ and a compactum $C\subseteq\partial\widehat M^n$ such that
$\widehat M^n-C$ is homeomorphic to $M^n$.  In that case $\widehat M^n$ is
called a \emph{manifold completion} of $M^n$.  Completable manifolds are also
called missing-boundary manifolds.  The analogous definitions in the PL and
DIFF categories use PL homeomorphisms and diffeomorphisms, respectively.

The first systematic classification of completable open manifolds of
dimension at least six that are simply connected at infinity was obtained by
Browder, Levine, and Livesay \cite{BLL65}.  Their proof helped make handle
theory the principal tool for manifold completion.  A full characterization,
allowing noncompact boundary and arbitrary end spaces, was completed by Gu
and Guilbault \cite{GG20}, following work of Siebenmann, Husch--Price,
Tucker, and O'Brien \cite{Si65,HP70,Tuc74,OB83}.  That theorem was the
culmination of more than five decades of work on the missing-boundary
problem.

\begin{theorem}[Manifold Completion Theorem \cite{GG20}]
\label{thm:completion}
An $n$-manifold $M^n$, where $n\ne4,5$, is completable if and only if
\begin{enumerate}[(1)]
\item\label{cond:tame} $M$ is inward tame;
\item\label{cond:stable} $M$ is peripherally $\pi_1$-stable at infinity;
\item\label{cond:sigma}
$\sigma_\infty(M)\in
\varprojlim\{\widetilde K_0(\mathbb Z[\pi_1(N)])\}$ is zero; and
\item\label{cond:tau}
$\tau_\infty(M)\in
\varprojlim{}^{1}\{\Wh(\pi_1(N))\}$ is zero.
\end{enumerate}
Here $N$ ranges over a cofinal sequence of clean neighborhoods of infinity.
\end{theorem}

For surfaces the corresponding statement is contained in their
classification \cite{Ke23,Ri63}.  For $3$-manifolds see
\cite{HP70,Tuc74,Ka87}.  Tucker's argument relies heavily on advances in
$3$-manifold theory during the 1970s \cite{GHM72,Ja71,Sc73}.  The conclusion
of \cref{thm:completion} fails in dimension $4$; counterexamples were found
independently by Weinberger and by Kwasik--Schultz
\cite{Wei87,KS88}.  In the TOP category, the dimension-$5$ conclusion holds
when the relevant stable peripheral groups are good in the sense of
Freedman--Quinn \cite{FQ90,GG20}.

Condition \ref{cond:stable} is too rigid for many examples arising in
geometric topology.  Davis's exotic universal covers, for example, are not
$\pi_1$-stable at infinity \cite{Dav83,Dav08}.  This led to the homotopy
analogue of manifold completion studied in \cite{Gui00,GT03,GT06,Gu20}.

\begin{definition}\label{def:pseudocollar}
A manifold neighborhood of infinity $N$ in a manifold $M$ is a
\emph{homotopy collar} if $\Fr_MN\hookrightarrow N$ is a homotopy
equivalence.  A \emph{pseudo-collar} is a homotopy collar that contains
arbitrarily small homotopy collar neighborhoods of infinity.  Finally, a
manifold is \emph{pseudo-collarable} if it contains a pseudo-collar
neighborhood of infinity.
\end{definition}

\begin{theorem}[Pseudo-collarability characterization \cite{Gu20}]
\label{thm:pseudocollar}
An $n$-manifold $M^n$, where $n\ge6$, is pseudo-collarable if and only if
\begin{enumerate}[(a)]
\item $M$ is inward tame;
\item $M$ is peripherally perfectly $\pi_1$-semistable at infinity; and
\item $\sigma_\infty(M)=0$.
\end{enumerate}
The TOP version also holds in dimension $5$ when the groups at infinity are
good.
\end{theorem}

The universal covers of all presently standard families of closed
aspherical manifolds are pseudo-collarable.  This motivates the following
question from \cite[Section~5]{Gui16}.

\begin{question}\label{q:pseudocollar-universal-cover}
Does every closed aspherical $n$-manifold, $n\ge4$, have a
pseudo-collarable universal cover?
\end{question}

The compactification problem considered here is weaker.  A closed subset
$Z$ of a metric ANR $\overline X$ is a \emph{$\Z$-set} if, for every open
set $U\subseteq\overline X$, the inclusion $U-Z\hookrightarrow U$ is a
homotopy equivalence.  A compactification $\overline X=X\sqcup Z$ is a
\emph{$\Z$-compactification} if $Z$ is a $\Z$-set in $\overline X$.
In this case $Z=\overline X-X$ is called the \emph{$\Z$-boundary} of
$X$.\footnote{Here ``$\Z$-boundary'' simply means the added $\Z$-set
in a $\Z$-compactification. No group action, cocompactness assumption,
or $\Z$-structure is part of this terminology.} The central open problem is the following.

\begin{question}\label{q:main}
Are Conditions \ref{cond:tame}, \ref{cond:sigma}, and \ref{cond:tau} of
\cref{thm:completion} sufficient for an open manifold to be
$\Z$-compactifiable?
\end{question}

The relation with pseudo-collarability leads to a more specific question.

\begin{question}\label{q:boundaryless-separation}
Does there exist a (boundaryless) open manifold that is
$\Z$-compactifiable but not pseudo-collarable?
\end{question}

The question has an affirmative answer in every classical
high-dimensional PL dimension.

\begin{theorem}\label{thm:boundaryless-example}
For every $n\ge6$ there exists a connected one-ended open PL
$n$-manifold $M^n$ which admits a finite-dimensional $\Z$-compactification
\[
\overline M=M\sqcup Z
\]
with $\overline M$ a compact ANR, but $M$ is not pseudo-collarable.
Moreover, the $\Z$-boundary of every $\Z$-compactification of $M$ is
non-ANR.
\end{theorem}

The examples are obtained from the reverse-plus manifolds of
Guilbault--Tinsley \cite{GT03}.  A relative based handle-slide and cancellation argument produces complementary
$2/3$-handle pairs with compatible markings on the two faces of each
block.  Relative cell-like lifting then supplies deformation retractions
whose tracks tend to zero in every retained inverse coordinate.  The
resulting inverse limit is the required compact ANR.  The restriction
$n\ge6$ is precisely the range in which the complementary-handle argument
applies; no assertion is made here in dimensions four or five.

We also record restrictions that explain the necessary wildness in the
examples of \cref{thm:boundaryless-example}.  An ANR $\Z$-boundary forces
stable fundamental group at each end, while a clean track-faithful
compact-expansion exhaustion forces pseudo-collarability.  These statements
are collected in \cref{prop:anr-expansion-restrictions}.

Question \ref{q:main} originated in the Hilbert cube setting
\cite{CS76,Ch76,Ge79,Wes90}.  It is false for arbitrary low-dimensional
polyhedra \cite{Gui01}, whereas Ferry proved a stable positive result for
finite-dimensional polyhedra \cite{Fe00}.  For manifolds Gu and Guilbault
proved that if an open $n$-manifold $M$, $n\ge5$, satisfies Conditions
\ref{cond:tame}, \ref{cond:sigma}, and \ref{cond:tau}, then
$M\times[0,1]$ is $\Z$-compactifiable; in fact it is completable
\cite[Theorem~1.2]{GG20}.  We use
\[
I=[-1,1]
\]
throughout.  The choice of length has no mathematical significance.

One of the main goals of this paper is to investigate when the
$\Z$-compactifiability of $M\times I$ implies the
$\Z$-compactifiability of $M$.  For a completion
\[
W=\widehat{M\times I}=(M\times I)\sqcup Z,
\]
put the central slice
\[
X=\overline{M\times\{0\}}^{\,W},
\qquad
Z_0=X-(M\times\{0\})=X\cap Z.
\]
The phrase ``$Z$ splits as $Z_0\times I$'' is not, by itself, enough to
control how the central slice approaches $Z$.  We refer the given identification $W-Z=M\times I$ as the \emph{marked product};
its distinguished intervals $\{m\}\times I$ are the \emph{marked
vertical fibers}, and the associated map
\[
 (m,t)\longmapsto(m,0)
\]
is the \emph{marked vertical projection}.  This terminology is used
only to distinguish the given product structure on $W-Z$ from an
abstract product decomposition of the added set $Z$.

The following theorem uses
the ambient compatibility that is needed in the proof.

\begin{theorem}\label{thm:main}
Suppose that $M^n$ is an open $n$-manifold, $n\ge5$, satisfying Conditions
\ref{cond:tame}, \ref{cond:sigma}, and \ref{cond:tau}.  Then $M\times I$
admits a manifold completion $W=(M\times I)\sqcup Z$.

Assume that this completion has a product-compatible splitting: there is a
homeomorphism
\[
\theta:Z_0\times I\longrightarrow Z
\]
such that
\[
\theta(z,0)=z \qquad (z\in Z_0),
\]
and
\[
\theta\bigl(Z_0\times(-1,1)\bigr)=\Int_{\partial W}Z.
\]
Assume, moreover, that the map
\begin{equation}\label{eq:compatible-projection}
q:W\longrightarrow X,
\qquad
q(m,t)=(m,0),\qquad
q\bigl(\theta(z,t)\bigr)=z,
\end{equation}
for $m\in M$, $z\in Z_0$, and $t\in I$, is continuous.
Then $X=(M\times\{0\})\sqcup Z_0$ is a $\Z$-compactification of
$M\times\{0\}\cong M$, and $Z_0$ is collared in $X$.  More precisely,
there is an embedding
\[
\lambda:Z_0\times[0,1)\longrightarrow X
\]
onto a neighborhood of $Z_0$ in $X$ with $\lambda(z,0)=z$, $z\in Z_0$.

For $n=4$, assume instead that $M\times I$ is supplied with a manifold
completion $W$ satisfying the same product-compatibility hypotheses.  For
example, the existence of $W$ follows from the dimension-$5$ TOP completion
theorem when $M\times I$ is inward tame, is peripherally $\pi_1$-stable with
good stable groups in the sense of Freedman--Quinn, and has vanishing
$\sigma_\infty$ and $\tau_\infty$.
\end{theorem}

Roughly speaking, the compatible splitting in \cref{thm:main} extends the
vertical projection to a retraction of the stabilized completion onto the
central slice.  The $\Z$-set homotopy then gives the required
$\Z$-compactification, and the normal-bundle argument gives the collar.  The
word ``compatible'' is essential: an abstract product homeomorphism of the
added set does not by itself determine how the marked middle slice approaches
the boundary. 

The $\Z$-compactification part requires less than a product-compatible
splitting.  The precise destabilization criterion is as follows.

\begin{theorem}\label{thm:retraction-intro}
Let $W=(M\times I)\sqcup Z$ be a compact metric ANR in which $Z$ is a
$\Z$-set, and define $X$ and $Z_0$ as above.  If there is an open
neighborhood $U$ of $Z_0$ in $W$ and a retraction
\[
r:U\longrightarrow U\cap X
\quad\text{with}\quad
r(U-Z)\subseteq M\times\{0\},
\]
then $X=(M\times\{0\})\sqcup Z_0$ is a $\Z$-compactification of $M$.
\end{theorem}

 In particular, a global extension of the vertical projection is sufficient,
but it is not necessary.  A second result below shows that, already for
$n\ge4$, the weaker abstract split $Z\cong A\times I$ forces $A$ to be an
ANR and forces stability at infinity.  When $n\ge6$, inward tameness and
$\sigma_\infty(M)=0$ then give an ordinary collar on $M$.  That collar need
not have boundary $A$ and does not
identify the topology inherited by the marked middle slice.

Neither \cref{thm:main} nor \cref{thm:retraction-intro} asserts that an
arbitrary stabilized completion has the required compatibility.  Producing
such compatibility is exactly the unresolved part of Question~\ref{q:main}.

The rest of the paper is organized as follows.  Section~2 reviews neighborhoods of
infinity, the two Chapman--Siebenmann obstructions, $\Z$-sets, and the
normal-bundle theorem used below.  Section~3 proves
\cref{thm:main,thm:retraction-intro}, gives metric, controlled-product, and
cell-like destabilization criteria, proves the abstract-split stability
result, and briefly revisits the Kwasik--Schultz examples.  Section~4
constructs the boundaryless examples of
\cref{thm:boundaryless-example}.  Section~5 separates ordinary
$\Z$-compactifiability from the controlled/nullity condition needed for
boundary swapping and applications to aspherical manifolds.

\subsection*{Acknowledgements}
I would like to thank Craig Guilbault for pointing out the error in
Theorem \ref{thm:main} of the initial version of the manuscript.  I also thank the
anonymous referee for a careful reading and for suggestions that substantially
improved the organization and accuracy of the paper.

\section{Preliminaries}\label{sec:preliminaries}

\subsection{Clean neighborhoods and ends}

Unless stated otherwise, manifolds may have boundary.  For $A\subseteq X$,
the point-set interior of $A$ in $X$ is denoted by $\Int_XA$, and
\[
\Fr_XA:=\overline A^{\,X}\cap\overline{X-A}^{\,X}
=\overline A^{\,X}-\Int_XA
\]
denotes its frontier.
A closed manifold is compact and boundaryless, whereas an open manifold is
noncompact and boundaryless.

A codimension-zero submanifold $Q\subseteq M$ is \emph{clean} if it is
closed in $M$ and $\Fr_MQ$ is a properly embedded, locally flat, and hence
bicollared codimension-one submanifold.  A clean neighborhood of infinity is
a clean codimension-zero submanifold $N\subseteq M$ for which
$\overline{M-N}$ is compact.  A nested sequence
\[
N_1\supseteq N_2\supseteq N_3\supseteq\cdots
\]
is \emph{cofinal} if every neighborhood of infinity contains some $N_i$.
After the usual enlargements of the complementary compacta, each component
of $N_i$ may be assumed noncompact and to have connected frontier; such an
$N_i$ is called a clean $0$-neighborhood of infinity.

An end $\varepsilon$ of $M$ is represented by a nested sequence of
components $N_1^{\varepsilon}\supseteq N_2^{\varepsilon}\supseteq\cdots$.
Choose a proper base ray $r$ meeting these components cofinally.  Inclusions,
together with change-of-basepoint paths along $r$, give
\begin{equation}\label{eq:propi}
\pi_1(N_1^{\varepsilon},r(1))
\longleftarrow
\pi_1(N_2^{\varepsilon},r(2))
\longleftarrow
\pi_1(N_3^{\varepsilon},r(3))
\longleftarrow\cdots .
\end{equation}
Its pro-isomorphism class is the fundamental pro-group
$\propi(\varepsilon,r)$.  Standard changes of cofinal sequence do not alter
that class; see \cite{Ge08,HR96}.

\subsection{Pro-groups at infinity}

Two inverse sequences are \emph{pro-isomorphic} if, after passage to
subsequences, they occur as the rows of a commuting ladder with maps in both
directions.  An inverse sequence is \emph{semistable} if it is pro-isomorphic
to one whose bonding maps are surjective.  It is \emph{pro-monomorphic} if it
is pro-isomorphic to one whose bonding maps are injective.  It is
\emph{stable} if it is pro-isomorphic to a constant inverse sequence
\[
H\xleftarrow{\id}H\xleftarrow{\id}H\xleftarrow{\id}\cdots .
\]
Equivalently, an inverse sequence is stable if it is both semistable and
pro-monomorphic.  An end is stable if one, hence every, representative of
its fundamental pro-group is stable.

An inverse sequence is \emph{perfectly semistable} if it is pro-isomorphic
to a sequence
\[
G_0\xleftarrow{\lambda_1}G_1\xleftarrow{\lambda_2}G_2
\xleftarrow{\lambda_3}\cdots
\]
of finitely generated groups and surjections whose kernels are perfect.
The peripheral versions apply these definitions to the appropriate
components of clean neighborhoods together with their boundary data; see
\cite{GG20,Gu20}.

\subsection{Tameness and the Wall obstruction}

A space $P$ is \emph{finitely dominated} if there are a finite polyhedron
$K$ and maps $u:P\to K$ and $d:K\to P$ such that
$d\circ u\simeq\id_P$.  It has \emph{finite homotopy type} if it is homotopy
equivalent to a finite polyhedron.  A locally finite polyhedron, or more
generally a locally compact ANR, is \emph{inward tame} if it has arbitrarily
small closed neighborhoods of infinity that are finitely dominated.

Let $\Lambda$ be a unital ring; in the applications below
$\Lambda=\mathbb Z[G]$.  The reduced projective class group
$\widetilde K_0(\Lambda)$ is the Grothendieck group of finitely generated
projective left $\Lambda$-modules modulo the subgroup generated by the class
of the free module $\Lambda$.  Wall associated to every connected finitely
dominated space $P$ an obstruction
\[
\sigma(P)\in\widetilde K_0(\mathbb Z[\pi_1(P)])
\]
that vanishes precisely when $P$ has finite homotopy type \cite{Wal65}.

For a space $P$ with finitely many path components
$P^1,\dots,P^k$, set
\[
\widetilde K_0\bigl(\mathbb Z[\pi_1(P)]\bigr)
:=\bigoplus_{j=1}^k
\widetilde K_0\bigl(\mathbb Z[\pi_1(P^j)]\bigr),
\qquad
\sigma(P):=(\sigma(P^1),\dots,\sigma(P^k)).
\]
If $P$ is inward tame, choose a nested cofinal sequence $\{N_i\}$ of
closed neighborhoods of infinity, each finitely dominated, and define
\[
\sigma_\infty(P)
:=(\sigma(N_1),\sigma(N_2),\dots)
\in
\varprojlim_i\widetilde K_0\bigl(\mathbb Z[\pi_1(N_i)]\bigr).
\]
The bonding maps are induced by inclusion, and the sum theorem for Wall's
obstruction gives compatibility.  The inverse limit is a subgroup of the
direct product, so its zero element is the coordinatewise-zero tuple.
Consequently, $\sigma_\infty(P)=0$ if and only if every $\sigma(N_i)$ is
zero.  By Wall's theorem and another application of the sum theorem, this is
equivalent to every closed polyhedral neighborhood of infinity having finite
homotopy type.

\subsection{The derived-limit torsion obstruction}

For an inverse sequence of abelian groups
\[
G_0\xleftarrow{\lambda_1}G_1\xleftarrow{\lambda_2}G_2
\xleftarrow{\lambda_3}\cdots,
\]
the derived limit is
\[
\varprojlim{}^1\{G_i,\lambda_i\}
:=
\frac{\prod_{i=0}^{\infty}G_i}
{\{(g_0-\lambda_1g_1,g_1-\lambda_2g_2,\dots):g_i\in G_i\}}.
\]
Pro-isomorphic inverse sequences of abelian groups have isomorphic derived
limits.

Suppose that $M$ contains a cofinal sequence $\{N_i\}$ of clean
neighborhoods of infinity for which every inclusion
$\Fr_MN_i\hookrightarrow N_i$ is a homotopy equivalence.  Put
$W_i=\overline{N_i-N_{i+1}}$.  The inclusion
$\Fr_MN_i\hookrightarrow W_i$ has a Whitehead torsion
\[
\tau_i=\tau(W_i,\Fr_MN_i)
\in\Wh(\pi_1(\Fr_MN_i))\cong\Wh(\pi_1(N_i)).
\]
For a disconnected frontier, this means
\[
\Wh(\pi_1(\Fr_MN_i))
:=\bigoplus_j\Wh(\pi_1(\Fr_MN_i^j)).
\]
The inclusion-induced maps give an inverse sequence of Whitehead groups, and
Chapman--Siebenmann's obstruction is
\[
\tau_\infty(M):=[(\tau_1,\tau_2,\dots)]
\in\varprojlim{}^1_i\Wh(\pi_1(N_i)).
\]
For the general definition and its invariance, see \cite{CS76,GG20}.

\subsection{\texorpdfstring{$\Z$}{Z}-sets, local fundamental groups, and normal bundles}

We record several facts used in the proof.  For a closed subset $Z$ of a
compact metric ANR $\overline X$, the following are equivalent
\cite{He75}:
\begin{enumerate}[(1)]
\item $Z$ is a $\Z$-set in $\overline X$;
\item there is a homotopy $H:\overline X\times[0,1]\to\overline X$ with
$H_0=\id$ and $H_t(\overline X)\cap Z=\varnothing$ for every $t>0$;
\item for every $\epsilon>0$ there is a self-map of $\overline X$ that is
$\epsilon$-close to the identity and whose image misses $Z$; and
\item $\overline X-Z$ is dense in $\overline X$, and $Z$ is $k$-LCC in
$\overline X$ for every $k\ge0$.
\end{enumerate}
Here $k$-LCC means that, for every $z\in Z$ and neighborhood $U$ of $z$,
there is a neighborhood $V\subseteq U$ of $z$ such that every map
$\partial I^{k+1}\to V-Z$ extends to a map $I^{k+1}\to U-Z$.

Using the preceding conventions for inverse systems and changes of
basepoint, let $N$ be a closed subset of a locally compact ANR $Y$ and
let $x\in N$.  The \emph{local fundamental pro-group of $Y-N$ at $x$}
is represented by
\[
   \{\pi_1(U-N,u_U)\}_{U\ni x},
\]
where $U$ ranges over sufficiently small connected neighborhoods of $x$,
and the basepoints $u_U\in U-N$ are chosen coherently along a path in
$Y-N$ approaching $x$, within a chosen local component of $Y-N$ at
$x$.  We say that this local fundamental pro-group is \emph{trivial},
respectively \emph{infinite cyclic}, if it is pro-isomorphic to the
constant inverse system $1$, respectively $\mathbb Z$.

We use the following normal-bundle theorem in the form relevant here.

\begin{theorem}[Codimension-one normal-bundle theorem]
\label{thm:normal-bundle}
Let $N$ be a closed subset of an ANR homology $n$-manifold $Y$ without
boundary.  Assume that $Y-N$ is a manifold, $N\times\mathbb R$ is an
$n$-manifold, and the following componentwise local condition holds: for
every $x\in N$ and every neighborhood $U$ of $x$, there is a neighborhood
$V\subseteq U$ of $x$ such that every loop in each component of $V-N$
contracts in the corresponding component of $U-N$.  Then $N$ has a normal
line bundle in $Y$.  For $n=4$ the statement uses the disk-embedding
hypotheses of Freedman--Quinn.
\end{theorem}

For $n\ge5$ this is the codimension-one case of
\cite[Section~3.4]{Qu79}; the dimension-$4$ formulation is
\cite[Theorem~9.3A, pp.~137--141]{FQ90}.  The hypothesis that $Y$ has no
boundary is essential. The following standard Davis example also fixes a useful point of notation.
Choose the one-ended arrangement
\[
Q^n=C_0^n\mathbin{\natural_\partial}(-C_1^n)
\mathbin{\natural_\partial}C_2^n
\mathbin{\natural_\partial}\cdots
\]
described in \cite[Figure~3]{Gui16}.  An arbitrary infinite boundary
connected sum is not well-defined without specifying its arrangement; in
particular, one-ended and two-ended arrangements differ.  The Davis manifold
is $D^n\cong\operatorname{int}Q^n$.  By
\cite[Example~31]{Gui16},
\[
\widehat Q^n=Q^n\sqcup\{\infty\}
\]
is a $\Z$-compactification, and the induced $\Z$-boundary for $D^n$ is
\[
Z_D=\partial Q^n\sqcup\{\infty\}.
\]
As explained in \cite[Exercise~10.3]{Gui16}, $Z_D$ is a non-ANR homology
$(n-1)$-manifold.  This illustrates why a product with $\mathbb R$ being a
manifold does not, without the remaining hypotheses of
\cref{thm:normal-bundle}, automatically produce a collar.

\section{Destabilizing a \texorpdfstring{$\Z$}{Z}-compactification}
\label{sec:destabilization}

\subsection{The central-slice retraction}

We first prove the splitting-free criterion stated in the introduction.
It separates the elementary $\Z$-set argument from the normal-bundle
argument needed for a collar.

\begin{proof}[Proof of \cref{thm:retraction-intro}]
The slice $M\times\{0\}$ is closed in $W-Z=M\times I$, so
$X-Z_0=M\times\{0\}$.  The open set $U$ is an ANR and $U\cap X$ is its
retract, hence an ANR.  Away from $Z_0$, the space $X$ is locally the
manifold $M$.  Thus $X$ is locally an ANR, and Hanner's local-to-global
theorem makes the compact metrizable space $X$ an ANR.

Let
\[
H:W\times[0,1]\longrightarrow W
\]
be a $\Z$-set homotopy for $Z$, so $H_0=\id_W$ and
$H_t(W)\subseteq W-Z=M\times I$ for every $t>0$.  Choose an open
$O\subseteq X$ with
\[
             Z_0\subseteq O,
             \qquad \overline O^{\,X}\subseteq U\cap X,
\]
and a continuous function $\lambda:X\to[0,1]$ which is positive on $Z_0$
and whose support $F$ is contained in $O$.  Compactness gives $\tau>0$ such
that $H(F\times[0,\tau])\subseteq U$.  After replacing $\lambda$ by
$\tau\lambda$, define
\[
K:X\times[0,1]\longrightarrow X,
\qquad
K(x,t)=
\begin{cases}
r(H(x,t\lambda(x))),&x\in F,\\
x,&x\notin F.
\end{cases}
\]
The formulas agree where $\lambda=0$, since there $r(x)=x$.  We have
$K_0=\id_X$.  If $t>0$ and $\lambda(x)>0$, then
$H(x,t\lambda(x))\in W-Z$ and therefore
$K(x,t)\in M\times\{0\}$.  If $\lambda(x)=0$, then $x\notin Z_0$.
Consequently $K_t(X)\cap Z_0=\varnothing$ for every $t>0$.  Thus $Z_0$ is
a $\Z$-set in $X$, and $X-Z_0\cong M$.
\end{proof}

The extension hypothesis has a useful metric sufficient condition.  If
$K\subseteq M$ is compact, write $M-K$ for its complement.

\begin{proposition}[Shrinking-fiber criterion]\label{prop:shrinking-fibers}
Let $W=(M\times I)\sqcup Z$ be a compact metric ANR $\Z$-compactification,
with compatible metric $\overline d$.  Suppose
\begin{equation}\label{eq:shrinking-fibers}
\lim_{K\nearrow M}
\sup\bigl\{\overline d((m,t),(m,0)):
m\in M-K,\ t\in I\bigr\}=0,
\end{equation}
where the limit is over compact subsets of $M$.  Then $Z=Z_0$, the vertical
projection extends to a retraction $R:W\to X$ by $R|_Z=\id_Z$, and $X$ is a
$\Z$-compactification of $M$.
\end{proposition}

\begin{proof}
Since $M\times I$ is dense in $W$, every $z\in Z$ is the limit of a sequence
$(m_i,t_i)$.  The points $m_i$ leave every compact subset of $M$; otherwise
a subsequence of $(m_i,t_i)$ would converge in $M\times I$.  By
\eqref{eq:shrinking-fibers},
\[
\overline d((m_i,t_i),(m_i,0))\longrightarrow0.
\]
Hence $(m_i,0)\to z$, so $z\in X$ and therefore $Z=Z_0$.  Define
\[
R(m,t)=(m,0),\qquad R(z)=z.
\]
For an arbitrary sequence $y_i\to z\in Z$, pass to the subsequences of terms
in $Z$ and in $M\times I$.  On the first, $R(y_i)=y_i\to z$; on the second,
the same shrinking estimate gives $R(y_i)\to z$.  Thus $R$ is sequentially
continuous at $Z$, and it is plainly continuous on $M\times I$.  All spaces
are metrizable, so $R$ is continuous.  It is a retraction and satisfies the
hypotheses of \cref{thm:retraction-intro}.
\end{proof}

The preceding criterion does not require the product lines supplied by a
controlled theorem to be the original vertical fibers.  What is required is
control that separates points of the added set and tends to zero relative to a
radial coordinate.

\begin{proposition}[Controlled-product projection]
\label{prop:controlled-product-projection}
Let $\overline{\mathcal H}$ be a compact metrizable space containing
$A\times I$ as a closed subset, and put
\[
   \mathcal H=\overline{\mathcal H}-(A\times I).
\]
Suppose there is a subspace $B_0\subset\mathcal H$ and a homeomorphism
\[
   \Phi:B_0\times I\xrightarrow{\cong}\mathcal H,
   \qquad \Phi(b,0)=b,
\]
such that
\[
   \overline{B_0}^{\,\overline{\mathcal H}}
   =B_0\sqcup(A\times\{0\}).
\]
Let $P$ be a metric space, let $j:A\hookrightarrow P$ be an embedding, and
let
\[
 q=(u,\rho):\overline{\mathcal H}\longrightarrow
 P\times[0,\infty)
\]
be continuous, where the target has the maximum metric and
\[
 q(a,t)=(j(a),0),\qquad
 \rho^{-1}(0)=A\times I,
 \qquad
 q^{-1}(j(a),0)\cap\overline{B_0}=\{(a,0)\}.
\]
If there is a constant $c<1$ such that
\begin{equation}\label{eq:conical-controlled-product}
 d\bigl(q(\Phi(b,s)),q(b)\bigr)<c\rho(b)
 \qquad(b\in B_0,\ s\in I),
\end{equation}
then projection along the $\Phi$-product lines extends continuously by
\[
 \overline r_\Phi(\Phi(b,s))=b,
 \qquad
 \overline r_\Phi(a,t)=(a,0).
\]
In particular, if $\overline{\mathcal H}$ is a neighborhood of the central
added set in a stabilized completion and
$B_0\sqcup(A\times\{0\})$ is its intersection with the central closure,
then, after shrinking to an invariant open neighborhood, this extension gives
the local retraction required in \cref{thm:retraction-intro}.
\end{proposition}

\begin{proof}
Let $y_k=\Phi(b_k,s_k)\to(a,t)$.  The radial coordinate in
\eqref{eq:conical-controlled-product} gives
\[
              \rho(b_k)\le \frac{\rho(y_k)}{1-c}\longrightarrow0.
\]
The horizontal coordinate then gives
$d_P(u(b_k),u(y_k))\to0$, and hence
$q(b_k)\to(j(a),0)$.  Compactness of $\overline{B_0}$ and the stated
point-separation of its zero fiber imply $b_k\to(a,0)$.  Therefore
$\overline r_\Phi(y_k)\to(a,0)$, proving continuity at the added set.
For the final assertion, choose an open set $O$ in the ambient completion
with $A\times\{0\}\subseteq O\subseteq\overline{\mathcal H}$ and put
\[
 U=O\cap\overline r_\Phi^{-1}(O\cap X).
\]
Then $U$ is an open neighborhood of the central added set,
$\overline r_\Phi(U)=U\cap X$, and the restriction is the retraction required
by \cref{thm:retraction-intro}.
\end{proof}

For a finite-dimensional, possibly non-ANR compactum $A$, the control base
in this proposition can still be chosen to be polyhedral: embed
$A\hookrightarrow\mathbb R^{2\dim A+1}$, extend its coordinate functions
over $\overline{\mathcal H}$, and adjoin distance from $A\times I$ as
$\rho$.  Thus non-ANRness of the added set is not a formal obstruction to this
conditional projection argument.  It does, however, prevent the separate
conclusion that $A\times\mathbb R$ is a manifold, since the slice retraction
would then make $A$ an ANR.

There is a more flexible criterion that does not require the added set to
split as a product.  Recall that a decomposition of a compact metric space is
\emph{upper semicontinuous} if its quotient is Hausdorff; equivalently, the
associated equivalence relation is closed.  A compactum is
\emph{cell-like} if it has the shape of a point, and a proper surjection is
cell-like if each of its point inverses is cell-like.

\begin{theorem}[Cell-like destabilization criterion]
\label{thm:cell-like-destabilization}
Let $W=(M\times I)\sqcup Z$ be a compact finite-dimensional metric ANR
$\Z$-compactification.  Suppose $W$ has an upper-semicontinuous decomposition
$\mathcal G$ such that
\begin{enumerate}[(i)]
\item every decomposition element meeting $M\times I$ is exactly a vertical
interval $\{m\}\times I$;
\item every remaining decomposition element is a cell-like subset of $Z$;
and
\item the quotient $Y=W/\mathcal G$ is finite-dimensional.
\end{enumerate}
Then, with $A$ the image of $Z$, the quotient is a
$\Z$-compactification
\[
Y=M\sqcup A
\]
of $M$.
\end{theorem}

\begin{proof}
Let $q:W\to Y$ be the quotient map.  Every point inverse of $q$ is cell-like,
and $q$ is proper.  Since $W$ is an ENR and $Y$ is finite-dimensional and
metrizable, \cite[Corollary~3.3]{Lac69} implies that $Y$ is an ENR.  The set
$A=q(Z)$ is compact, $q^{-1}(A)=Z$, and
\[
Y-A\cong (M\times I)/(\{m\}\times I)\cong M.
\]

Let $U\subseteq Y$ be open.  The restrictions
\[
q^{-1}(U)\longrightarrow U
\quad\text{and}\quad
q^{-1}(U-A)\longrightarrow U-A
\]
are proper cell-like maps between ENRs, so
\cite[Theorem~1.2]{Lac69} says that both are proper homotopy equivalences.
On the other hand, the $\Z$-set property gives a homotopy equivalence
\[
q^{-1}(U)-Z\longrightarrow q^{-1}(U).
\]
Since $q^{-1}(U)-Z=q^{-1}(U-A)$, the commutative square formed by these
three maps shows that $U-A\hookrightarrow U$ is a homotopy equivalence.
This holds for every open $U$, so $A$ is a $\Z$-set in $Y$.
\end{proof}

The finite-dimensional quotient condition is automatic for the natural
vertical decomposition when that decomposition is upper semicontinuous.

\begin{corollary}\label{cor:usc-splitting}
Let $M^n$ be an open $n$-manifold and let
$W=(M^n\times I)\sqcup Z$ be a compact $(n+1)$-manifold completion, with
$Z\subseteq\partial W$.  Suppose there is a
homeomorphism $\theta:A\times I\to Z$ and let
$\mathcal G$ consist of the intervals $\{m\}\times I$, $m\in M$, together
with $\theta(\{a\}\times I)$, $a\in A$.  If $\mathcal G$ is upper
semicontinuous, then $W/\mathcal G=M\sqcup A$ is a
$\Z$-compactification of $M$.  If, in addition,
\[
\theta(A\times(-1,1))=\Int_{\partial W}Z,
\]
then $A$ is collared in $W/\mathcal G$ for $n\ge5$; the same holds for
$n=4$ under the Freedman--Quinn hypotheses.
\end{corollary}

\begin{proof}
The quotient is compact metrizable.  The embedded copy
$\theta(A\times\{0\})\subseteq\partial W$ meets every decomposition element
in $Z$ once.  Hence the restriction of the quotient map to this copy is a
continuous bijection onto $q(Z)$; it is a homeomorphism because its domain is
compact and the quotient is Hausdorff.  In particular, $A$ has dimension at
most $n$.
Choose a compact exhaustion $K_1\subseteq K_2\subseteq\cdots$ of $M$.
In the quotient, $A$ and all $K_i$ are closed, and
\[
W/\mathcal G=A\cup\bigcup_iK_i.
\]
The countable closed-sum theorem for covering dimension
\cite[Theorem~4.3.7]{vM89} gives
$\dim(W/\mathcal G)\le n$.  Thus
\cref{thm:cell-like-destabilization} applies.

For the final assertion,
$A\times\mathbb R\cong\Int_{\partial W}Z$ is an $n$-manifold.  The doubling
argument in \cite[Corollary~4]{AG99} gives an ANR homology $n$-manifold without
boundary when $n\ge5$; in dimension $4$ use the pre-disjoint-disks portion
of that proof together with
\cite[Theorem~9.3A]{FQ90}.  The $\Z$-set
property gives trivial local fundamental pro-groups along $A$.  The
normal-bundle argument from the proof of \cref{thm:main} therefore gives a
two-sided normal line bundle and a collar on either side.
\end{proof}

\begin{proposition}
\label{prop:split-corona-stability}
Let $M^n$ be an inward-tame open manifold, $n\ge4$, and suppose a manifold
completion
\[
                  W=(M\times I)\sqcup Z
\]
admits a homeomorphism $\theta:A\times I\to Z$ satisfying
\[
       \theta(A\times(-1,1))=\Int_{\partial W}Z.
\]
No compatibility with the marked vertical projection is assumed.  Then $A$
is a compact ANR and every end of $M$ has stable fundamental group
at infinity.  If in addition $n\ge6$ and
$\sigma_\infty(M)=0$, every end is collared and $M$ is the interior of a
compact manifold.  The derived-limit obstruction $\tau_\infty(M)$ is not
needed for this conclusion.
\end{proposition}

\begin{proof}
The space $A\times(-1,1)$ is an open subset of the $n$-manifold
$\partial W$, hence is an ANR.  Retraction to the zero slice makes $A$ an
ANR; it is compact because it is a closed slice of the compactum $Z$.
Consequently $Z\cong A\times I$ is a compact ANR.  Since $Z$ is a closed
subset of the boundary of $W$, a boundary collar instantly pushes $W$ off
$Z$.  Thus $W$ is a $\Z$-compactification of $M\times I$.

Inward tameness of the open manifold $M$ gives finitely many ends and
semistable fundamental group at each end
\cite[Theorem~1.2]{GT03}.  Crossing with the compact interval preserves
inward tameness, the end set, and the corresponding pro-fundamental groups.
The shape-of-the-end theorem identifies each end shape of $M\times I$ with
the corresponding component $A_j\times I$ of its $\Z$-boundary
\cite[Theorem~8.4]{Gui16}.  Semistability permits the componentwise pointed
pro-fundamental-group comparison, including the change from the boundary
basepoint to a base ray \cite[Theorem~7.13(2)]{Gui16}.  Each
$A_j\times I$ is a compact ANR and therefore has a constant, hence stable,
shape pro-fundamental group.  Thus the pro-fundamental group of each end of
$M\times I$, and therefore of $M$, is stable.

When $n\ge6$, Siebenmann's Collaring Theorem applies componentwise: inward tameness,
stability, and $\sigma_\infty(M)=0$ give an open collar at every end
\cite{Si65}.  There are only finitely many ends, so adjoining their compact
manifold frontiers gives a compact manifold with interior $M$.  No
$\varprojlim^{1}\Wh$ obstruction occurs after the system is stable.
\end{proof}

\begin{remark}[Exact scope of the abstract split]
\label{rem:split-corona-scope}
In \cref{prop:split-corona-stability}, the space $A$ is only the
abstract first factor in a chosen product decomposition
\[
   \theta:A\times I\longrightarrow Z.
\]
It is not assumed that the zero slice $\theta(A\times\{0\})$ coincides
with the actual subset
\[
   Z_0=\overline{M\times\{0\}}^{\,W}\cap Z,
\]
nor is any compatibility with the distinguished vertical fibers of
$M\times I$ assumed.  Thus, even if $A$ happens to be homeomorphic to
$Z_0$, the abstract product structure does not determine the topology
inherited by the marked middle slice.

Under the additional hypotheses $n\ge6$ and
$\sigma_\infty(M)=0$, the proposition produces \emph{some} manifold
boundary $B$.  Both $A$ and $B$ represent the same end shape and, being
compact ANRs, are homotopy equivalent componentwise.  This does not
imply that
\[
   X=\overline{M\times\{0\}}^{\,W}
     =(M\times\{0\})\sqcup Z_0
\]
is a $\Z$-compactification of $M$, nor that the resulting collar has
boundary $A$.  The example below shows that the marked vertical
projection can fail to extend even when $Z$ admits an abstract product
decomposition.

There are two further limitations.  First, the final Siebenmann step is
high-dimensional; the Kwasik--Schultz \cite{KS88} examples show that inward
tameness, stability, and $\sigma_\infty=0$ alone do not imply
collarability in dimension four.  The ANR and stability conclusions do
remain valid there.  Second, the split hypothesis itself forces $A$ to
be an ANR.  It therefore cannot realize a construction whose intended
$\Z$-boundary is genuinely non-ANR.  Such a construction must use the
local retraction or cell-like criteria without an ambient product
$A\times(-1,1)$ in $\partial W$.
\end{remark}

\begin{example}[An abstract split need not be compatible]
\label{ex:winding-split}
Let $W=D^2\times I$, $Z=S^1\times I$, and $M=\operatorname{int}D^2$.
Choose a continuous function $g:[0,1)\to\mathbb R$ that is zero near $0$
and tends to infinity at $1$, and identify $M\times I$ with $W-Z$ by
\[
h(re^{i\vartheta},t)
=\bigl(re^{i(\vartheta+t g(r))},t\bigr).
\]
The middle slice is standard, so $X=D^2\times\{0\}$,
$Z_0=S^1\times\{0\}$, and abstractly $Z\cong Z_0\times I$.  Choose
$r_j\to1$ with $g(r_j)=2\pi j+\alpha$, where
$\alpha\notin2\pi\mathbb Z$.  The points
\[
h(r_je^{i\vartheta},0)
\quad\text{and}\quad
h(r_je^{i\vartheta},1)
\]
lie on the same vertical fiber in the given product parametrization, but
converge to $(e^{i\vartheta},0)$ and
$(e^{i(\vartheta+\alpha)},1)$, respectively.  Those points lie on different
prescribed boundary fibers.  Hence the vertical equivalence relation is not
closed.

The failure of projection extension can also be seen directly.  Fix
$b=(e^{i\beta},1)$ and choose distinct residues
$\alpha,\alpha'\pmod{2\pi}$.  Take $r_j,r'_j\to1$ with
$g(r_j)=2\pi j+\alpha$ and $g(r'_j)=2\pi j+\alpha'$.  Then
\[
h(r_je^{i(\beta-\alpha)},1)\longrightarrow b,
\qquad
h(r'_je^{i(\beta-\alpha')},1)\longrightarrow b,
\]
whereas their images under the vertical projection converge in the middle
slice to the distinct points $(e^{i(\beta-\alpha)},0)$ and
$(e^{i(\beta-\alpha')},0)$.  Thus the projection has no continuous extension,
despite the abstract product decomposition of $Z$.  
\end{example}

\subsection{A product-compatible split and a collar}

We now prove the stronger conclusion of \cref{thm:main}.  The proof makes
explicit each hypothesis of the normal-bundle theorem.

\begin{proof}[Proof of \cref{thm:main}]
Thus let $M^n$, $n\ge5$, be inward tame with
$\sigma_\infty(M)=0$ and $\tau_\infty(M)=0$, and let
$W=(M\times I)\sqcup Z$ be a completion.  Put
$X=\overline{M\times\{0\}}^{\,W}$ and $Z_0=X\cap Z$.  Assume that
$\theta:Z_0\times I\to Z$ is a homeomorphism satisfying
$\theta(z,0)=z$ and
$\theta(Z_0\times(-1,1))=\Int_{\partial W}Z$, and that the map
$q:W\to X$ defined in \eqref{eq:compatible-projection} is continuous.  In
dimension $4$, begin instead with the supplied completion and the same
product-compatibility assumptions.
For $n\ge5$, the existence of the completion
$W=\widehat{M\times I}=(M\times I)\sqcup Z$ follows from
\cite[Theorem~1.2]{GG20}.  In dimension $4$, $W$ is part of the hypothesis;
the dimension-$5$ good-group completion criterion stated in the theorem is
one way of constructing it.
Since $Z$ is a closed subset of the boundary of the manifold $W$, a boundary
collar pushes $W$ instantly off $Z$; thus $W$ is a $\Z$-compactification of
$M\times I$.

The map $q$ in \eqref{eq:compatible-projection} is a retraction onto $X$ and
sends $M\times I$ into $M\times\{0\}$.  Therefore
\cref{thm:retraction-intro} shows that
$X=(M\times\{0\})\sqcup Z_0$ is a $\Z$-compactification of $M$.

It remains to construct the collar.  Let
\[
D(X)=X\cup_{Z_0}X
\]
be the double.  For $n\ge5$, the doubling argument of
\cite[Corollary~4]{AG99} shows that it is an ANR homology $n$-manifold without
boundary.  In dimension $4$, the part of that proof preceding the
disjoint-disks argument, \cite[pp.~1270--1271]{AG99}, gives the same ANR
homology-manifold conclusion.  Furthermore,
\[
D(X)-Z_0\cong M\sqcup M
\]
is an $n$-manifold.  The compatible splitting also gives
\[
Z_0\times(-1,1)\cong\Int_{\partial W}Z,
\]
so $Z_0\times\mathbb R$ is an $n$-manifold.

Because $Z_0$ is a $\Z$-set in $X$, it is $1$-LCC there.  To verify the
componentwise local condition in the double, let $z\in Z_0$ and let $U$ be a
neighborhood of $z$ in $D(X)$.  In each copy of $X$, local contractibility
and the $1$-LCC property give a smaller neighborhood $V$ such that every loop
in $V-Z_0$ contracts in the corresponding component of $U-Z_0$.  Since
\[
D(X)-Z_0=M_+\sqcup M_-,
\]
each loop in the punctured smaller neighborhood lies entirely in one of
these labelled sides.  Hence the local fundamental pro-group on each side is
trivial.

The hypotheses of \cref{thm:normal-bundle} are now satisfied with $Y=D(X)$
and $N=Z_0$.  Hence $Z_0$ has a normal line bundle in $D(X)$.  The two
globally labelled components $M_+$ and $M_-$ orient that line bundle, so it
is trivial.  Its two half-bundles lie in the corresponding copies of $X$ and
give product collars $Z_0\times[0,1)$, as required.
\end{proof}


\subsection{Dimension-four candidates}

Kwasik and Schultz constructed dimension-$4$ counterexamples to the
manifold completion theorem.  Briefly, for $G=C_2$ and a suitable class
$\alpha\in L_4^h(G;1)$, with trivial orientation character, their
equivariant surgery construction gives a
compact cobordism $T(V,\alpha)$.  Concatenating copies in both directions
produces a free $G$-manifold $T_\infty(V,\alpha)$ properly homotopy
equivalent to $S^3\times\mathbb R$; its two-point compactification is $S^4$
with a semifree $G$-action \cite[Section~2]{KS88}; see also
\cite[Chapter~13A]{Wal99} for the surgery notation.  Put
\[
                    M_\alpha=T_\infty(V,\alpha)/C_2.
\]

A connected clean end neighborhood $U$ with connected frontier is a
\emph{$1$-neighborhood} if the inclusion induces an isomorphism
\[
\pi_1(\Fr U)\longrightarrow\pi_1(U).
\]

\begin{proposition}\label{prop:KS-candidates}
The open $4$-manifolds $M_\alpha$ are inward tame, have stable fundamental
group $C_2$ at both ends, and
satisfy $\sigma_\infty(M_\alpha)=\tau_\infty(M_\alpha)=0$, but are not
pseudo-collarable.
\end{proposition}

\begin{proof}
Kwasik--Schultz show that the ends are tame and stable and that both
projective end obstructions vanish \cite[Theorem~2.1(2) and p.~449]{KS88}.
Thus $\sigma_\infty=0$.  The stable end group is $C_2$, and
$\Wh(C_2)=0$, so the derived-limit torsion obstruction vanishes as well.
They also prove directly that neither end admits a
$1$-neighborhood \cite[p.~449]{KS88}; consequently neither end can contain a
homotopy collar, since every homotopy-collar neighborhood is a
$1$-neighborhood.
\end{proof}

By \cref{prop:shrinking-fibers}, any completion of
$M_\alpha\times I$ satisfying \eqref{eq:shrinking-fibers} would make
$M_\alpha$ $\Z$-compactifiable and hence would settle the
Kwasik--Schultz case positively.  Existing stabilized product and
absorption results do not provide that marked shrinking control.  Conversely,
the absence of $1$-neighborhoods does not obstruct a $\Z$-compactification
with a wild $\Z$-boundary.  Consequently the $M_\alpha$ remain candidates;
neither their $\Z$-compactifiability nor their
non-$\Z$-compactifiability is proved here.

\begin{remark}
After the first version of this paper appeared on arXiv, there have
been further developments related to the questions considered above.
In particular, in a companion paper \cite{GuKS}, the author proves
that the non-desuspendable Kwasik--Schultz open $4$-manifolds are
$\mathcal Z$-compactifiable.  
\end{remark}

\subsection{Pseudo-collar consequences}

Before constructing the high-dimensional examples, we record two general
consequences.  They explain why the $\Z$-boundary must be wild and why the most
direct compact-expansion construction cannot produce such an example.

\begin{proposition}[ANR $\Z$-boundaries and direct compact expansions]
\label{prop:anr-expansion-restrictions}
Let $M$ be a boundaryless open manifold.
\begin{enumerate}[(i)]
\item If $M\sqcup Z$ is a $\Z$-compactification and $Z$ is an ANR, then
 every end of $M$ has stable pro-fundamental group.  If $\dim M\ge6$ and
 $\sigma_\infty(M)=0$, every end is collared.
\item Suppose $M$ has a clean compact codimension-zero exhaustion
\[
 K_0\subset\Int K_1\subset\Int K_2\subset\cdots
\]
such that $K_{i+1}$ strongly deformation retracts onto $K_i$ by a homotopy
$H^i$ satisfying
\[
                         H^i_1H^i_t=H^i_1
                         \qquad(0\le t\le1).
\]
Then the neighborhoods $M-\Int K_i$ are homotopy collars and $M$ is
pseudo-collarable.
\end{enumerate}
\end{proposition}

\begin{proof}
A $\Z$-compactification makes $M$ inward tame.  The shape-of-the-end
comparison identifies the pro-fundamental group of each end with the pointed
shape pro-fundamental group of the corresponding component of $Z$
\cite[Theorems~7.13(2) and~8.4]{Gui16}; inward tameness supplies the needed
semistability \cite[Theorem~1.2]{GT03}.  A compact ANR has constant shape
pro-homotopy type.  This proves stability, and Siebenmann's Collaring Theorem
gives the last assertion in (i).

For (ii), concatenating the stage homotopies gives a track-faithful
retraction $K_j\searrow K_i$ \cite[Lemma~4.8]{GG25}.  It restricts to a
retraction
\[
                 K_j-\Int K_i\searrow\partial K_i;
\]
indeed a track entering $\Int K_i$ would, by track faithfulness, end at that
point, whereas the retraction fiber over an interior point is a singleton
\cite[Lemma~4.7(2)]{GG25}.  Every compact homotopy in
$M-\Int K_i$ lies in a finite shell.  Hence
$\partial K_i\hookrightarrow M-\Int K_i$ is a weak equivalence and, by
Whitehead's theorem for manifold ANRs, a homotopy equivalence.  These
homotopy collars are cofinal.
\end{proof}

The preceding proposition shows that any boundaryless
$\Z$-compactifiable manifold which is not pseudo-collarable must have
a genuinely wild $\Z$-boundary: in particular, its $\Z$-boundary
cannot be an ANR.  The next section constructs such examples from the
reverse-plus manifolds of Guilbault--Tinsley.

\section{Boundaryless examples}\label{sec:boundaryless-examples}

We now prove \cref{thm:boundaryless-example}.  The construction uses the
reverse-plus manifolds of Guilbault--Tinsley, but the ordinary deformation
retractions in their proof must be chosen with additional control on both
faces of every block.  We first record the inverse-limit form of that
control.

For $n\ge6$, write the Guilbault--Tinsley manifold \cite{GT03} as
\[
 M_*^n=(S^{n-2}\times B^2)\cup A_0\cup A_1\cup\cdots,
 \qquad
 K_i=(S^{n-2}\times B^2)\cup A_0\cup\cdots\cup A_{i-1},
\]
so that $\partial K_i=\Gamma_i$.  Their construction identifies
\[
 \pi_1\Gamma_i=G_i=
 \left\langle a_0,\ldots,a_i\ \middle|\
 a_s=[a_s,a_{s-1}],\ 1\le s\le i\right\rangle
\]
and identifies passage across $A_i$ with the epimorphism
$G_{i+1}\twoheadrightarrow G_i$ that kills $a_{i+1}$.  For $i\ge1$,
borrow the $2$-handle from the preceding block and put
\[
 A_i'=A_i\cup k_{i-1}^2,
 \qquad
 F_i=\Gamma_i\cup k_{i-1}^2.
\]
After choosing collars, these are literal PL subspaces satisfying
\[
 K_i\cap A_i'=F_i,
 \qquad
 K_i\cup A_i'=K_{i+1},
 \qquad
 A_{i-1}'\cap A_i'=F_i.
\]
Let
\[
   J_i=F_i\cap F_{i+1}.
\]
We refer to $J_i$ as the \emph{seam}; geometrically, it is the
surviving coattaching part of the new handle $k_i^2$ in the old
frontier.

We use the word \emph{marking} in the following elementary sense.
A marking of a face of a block is a specified PL identification of that
face with the corresponding prescribed subpolyhedron $F_i$ or $F_{i+1}$,
together with the chosen handle characteristic coordinates.  When the
two faces meet along the protected coattaching subpolyhedron
$J_i=F_i\cap F_{i+1}$, the markings are called \emph{compatible} if
their restrictions to $J_i$ agree pointwise.  Accordingly, a map or
deformation is said to respect the markings when, after these fixed
identifications, its restrictions are the stated maps on the literal
faces, rather than merely being conjugate to them by unspecified
self-homeomorphisms.

Proposition~4.4 and Figures~4--5 of \cite{GT03} provide a strong
deformation retraction $A_i'\to F_i$ and a finite Figure~5 core
\[
 C_i=\Gamma_i'\cup b_{i-1}^{n-2}\subset F_i.
\]

The issue is to choose those retractions coherently on $F_i$, $F_{i+1}$,
and $J_i$.  The core choices below will satisfy
\begin{equation}\label{eq:gt-seam-in-cores}
 J_i\subset C_i\cap C_{i+1},
\end{equation}
so that pointwise fixation of $J_i$ by both adjacent retractions is
well-defined.

For retractions
\[
 \rho_i:K_{i+1}\longrightarrow K_i
\]
that fix $K_i$, write, for $j\le k$,
\[
 \rho_{j,k}=\rho_j\rho_{j+1}\cdots\rho_{k-1}:K_k\longrightarrow K_j,
 \qquad \rho_{j,j}=\id_{K_j},
\]
and put
\[
 L=\varprojlim(K_i,\rho_i).
\]
Give each $K_j$ a metric $\operatorname{dist}_j\le1$.

\begin{proposition}\label{prop:gt-retained-control}
Suppose that, for all sufficiently large $i$, the block $A_i'$ has a
retraction $R_i:A_i'\to F_i$ and a strong deformation retraction
\[
 H_i:A_i'\times I\longrightarrow A_i'
\]
from the identity to $R_i$, fixed on $F_i$.  Paste $H_i$ to the stationary
homotopy on $K_i$, and paste $R_i$ to $\id_{K_i}$.  Denote the resulting
terminal retraction on $K_{i+1}$ by $\rho_i$.  Suppose there are integers
$a_i\le i$ and numbers $\varepsilon_i>0$ such that
\[
 a_i\longrightarrow\infty,
 \qquad
 \varepsilon_i\longrightarrow0,
\]
and, whenever $j\le a_i$,
\begin{equation}\label{eq:gt-retained-coordinate}
 \operatorname{dist}_j\!
 \left(\rho_{j,i+1}H_i(x,t),\rho_{j,i+1}(x)\right)<\varepsilon_i
 \qquad (x\in K_{i+1},\ t\in I).
\end{equation}
Then $L$ is a finite-dimensional compact ANR, the canonical map
\[
 e:M_*^n=\bigcup_i K_i\longrightarrow L,
\]
which sends a point to its eventual inverse-limit thread, identifies $M_*^n$
with an open dense subset of $L$ carrying its manifold topology, and
$L-M_*^n$ is a $\Z$-set.  In particular $L=M_*^n\sqcup(L-M_*^n)$ is a
finite-dimensional $\Z$-compactification of $M_*^n$.
\end{proposition}

\begin{proof}
Discarding and reindexing a finite initial segment, we may assume that the
hypotheses hold for every $i$.  The system is fully insulated in the sense
of \cite[Theorem~1.7]{GG25}.  Indeed, if $x\in K_m$, choose $i\ge m+2$ and
an open neighborhood
\[
 U\subset \Int K_{i-1}\subset K_i-F_i
\]
of $x$.  Since $\rho_i(K_{i+1}-K_i)\subset F_i$, induction gives
$\rho_{i,k}^{-1}(U)=U$ for every $k\ge i$.  The cited theorem and
\cite[Proposition~1.11]{GG25} therefore show that the canonical map $e$
identifies the direct union with an open subspace of $L$, with its original
direct-limit, hence manifold, topology.

Let $e_i:K_i\to L$ be the canonical section, let $\pi_i:L\to K_i$ be
projection, and put $q_i=e_i\pi_i$.  The maps $q_i$ converge uniformly to
the identity in the weighted inverse-limit metric
\[
 D(x,y)=\sum_{j\ge1}2^{-j}\operatorname{dist}_j(x_j,y_j).
\]
The stage homotopy induces a homotopy from $q_{i+1}$ to $q_i$ by
\[
 \widehat H_i(x,t)=e_{i+1}\bigl(H_i(\pi_{i+1}(x),t)\bigr).
\]
By \eqref{eq:gt-retained-coordinate}, each of its tracks has diameter at
most
\[
 2\varepsilon_i+2^{-a_i}.
\]
For fixed $N$, concatenate the $\widehat H_i$, $i\ge N$, on consecutive
dyadic intervals accumulating at time zero, and give time zero the value
$\id_L$.  The displayed bound and the convergence of $q_i$ make the
concatenation continuous at time zero.  Each piece is compared with the
original inverse-limit thread, rather than with the preceding piece, so the
coordinate errors do not accumulate.  It is a homotopy from $\id_L$ to
$q_N$, and at every positive time its image lies in an eventual finite
stage $e_i(K_i)\subset M_*^n$.  Thus the added set $L-M_*^n$ is homotopy
negligible.  Since $M_*^n$ is an ANR, \cite[Proposition~5.13]{GG25} makes
$L$ an ANR; $L$ is compact and metrizable because it is an inverse limit of
compact metrizable spaces.  The same instant homotopy shows that $L-M_*^n$
is a $\Z$-set.  Finally, $\dim L\le n$ by the dimension inequality for
inverse limits of compact metric spaces.
\end{proof}

Following Cohen, a retraction of compact polyhedra is \emph{collapsible} if
all of its point inverses are collapsible polyhedra \cite[Section~8]{Coh69}.
Every collapsible retraction is cell-like.

\begin{theorem}[Marked Guilbault--Tinsley realization]
\label{thm:gt-marked-realization}
For every $n\ge6$, the choices in the Guilbault--Tinsley construction may be
made so that there are finite subpolyhedra $C_i\subset F_i$, core sections
$s_i:C_i\hookrightarrow F_i$, and onto collapsible retractions
\[
 d_i:F_i\longrightarrow C_i,
 \qquad
 c_i:A_i'\longrightarrow C_i
\]
such that, for suitable maps $\bar f_i:C_{i+1}\to C_i$,
\begin{equation}\label{eq:gt-marked-face-square}
 c_i|_{F_i}=d_i,
 \qquad
 c_i|_{F_{i+1}}=\bar f_i d_{i+1},
 \qquad
 d_i|_{J_i}=d_{i+1}|_{J_i}=\id_{J_i}.
\end{equation}
Consequently $M_*^n$ admits a finite-dimensional compact ANR
$\Z$-compactification.
\end{theorem}

The middle equality in \eqref{eq:gt-marked-face-square} is the literal
commutativity of the marked face square
\[
\begin{CD}
 F_{i+1} @>{d_{i+1}}>> C_{i+1}\\
 @VVV @VV{\bar f_i}V\\
 A_i' @>{c_i}>> C_i,
\end{CD}
\]
where the left vertical arrow is inclusion.  The proof uses standard
complementary-handle and collapse constructions from Rourke--Sanderson and
Cohen.

\begin{lemma}[Complementary-pair collapse]
\label{lem:gt-complementary-pair-collapse}
In the standard complementary-handle pair of \cite[Lemma~6.6]{RS72}, let
$1\le r\le n-1$, put $X=I^r$ and $Z=I^{n-r-1}$, and choose a triangulation
for which $X\searrow\{x_0\}$ with $x_0\in\Int X$.  Put
\[
 B=X\times[-1,3]\times Z,
 \qquad K=X\times[1,3]\times Z,
 \qquad L=X\times[-1,1]\times Z,
\]
\[
 P=X\times\partial([1,3]\times Z)
   \cup(\{x_0\}\times[1,3]\times Z),
\]
and
\[
 R=X\times\bigl(([-1,3]\times\partial Z)\cup(\{3\}\times Z)\bigr).
\]
For every prescribed collapsible retraction $d:K\to P$ that fixes $P$,
there are a collapsible retraction $\varrho:B\to R$ and a PL map
$e:P\to R$ such that
\[
 \varrho|_K=e\,d,
 \qquad
 \varrho|_R=\id_R.
\]
In particular, both retractions fix every subpolyhedron of $K\cap R$.
\end{lemma}

\begin{proof}
Cohen's product-collapse lemma \cite[Lemma~8.3]{Coh69} gives
\[
 K\searrow X\times\partial([1,3]\times Z)
 \cup(\{x_0\}\times[1,3]\times Z)=P,
\]
so prescribed retractions of the stated kind exist.  Fix one such $d$.
Since $K\cap L=X\times\{1\}\times Z\subset P$, paste $d$ to the identity
on $L$.  The pasted map is a collapsible retraction of $B$ onto $P\cup L$;
its nonsingleton point inverses are precisely the point inverses of $d$.
Collapse the latter space to $R$ in the following order:
\[
 X\times[-1,1]\times Z
 \searrow
 (X\times[-1,1]\times\partial Z)\cup(X\times\{1\}\times Z),
\]
\[
 X\times\{1\}\times Z
 \searrow
 (X\times\{1\}\times\partial Z)\cup(\{x_0\}\times\{1\}\times Z),
\]
\[
 \{x_0\}\times[1,3]\times Z
 \searrow
 \{x_0\}\times\bigl(([1,3]\times\partial Z)\cup(\{3\}\times Z)\bigr).
\]
The first collapse is relative to $K\cap L$; at the second and third stages
the rest of the polyhedron meets the moving face in the displayed target.
Hence the sequences extend over the union and are relative to $R$.  Let
$e$ be the restriction to $P$ of the last three collapses.  Cohen's
Theorem~8.1 and Lemma~8.6 show that the composite $\varrho$ is collapsible
and give the literal factorization $\varrho|_K=e\,d$.
\end{proof}

\begin{lemma}\label{lem:gt-one-sided-prism}
Let $E\searrow E_0$ be a specified finite sequence of elementary collapses
with chosen standard ball-face retractions, and let $a:E\to E_0$ be their
composite.  There is a collapsible retraction
\[
 \widehat a:E\times I\longrightarrow (E\times\{1\})\cup(E_0\times I)
\]
whose restriction to $E\times\{0\}$ is $x\mapsto(a(x),0)$.
\end{lemma}

\begin{proof}
It is enough to treat one elementary collapse $E=E_0\cup B^m$, where
$F^{m-1}=E_0\cap B$ is a ball face.  Write $B=F\times I_s$.  Parametrize
the left edge followed by the top edge of the normal $(s,t)$-square by
$[0,2]$ and set
\[
 f(s,t)=
 \begin{cases}
  2t,& t\le s,\\
  s+t,& t\ge s.
 \end{cases}
\]
After subdivision along $t=s$ and $f^{-1}(1)$ this is PL.  It fixes the
left and top edges, sends the bottom edge to its lower-left endpoint, and
has polygonal-arc point inverses.  Its product with $F$ is therefore a
collapsible retraction of $B\times I$ onto
$(B\times\{1\})\cup(F\times I)$ with the prescribed bottom restriction.
Perform these prism retractions in the chosen collapse order and apply
\cite[Lemma~8.6]{Coh69}.
\end{proof}

\begin{lemma}[Coattaching-relative Figure~5 collapse]
\label{lem:gt-coattaching-figure5}
For a PL $2$-handle $D^2\times D^{n-2}$, the Figure~5 collapse from the
attaching side onto
\[
 (D^2\times S^{n-3})\cup(\{0\}\times D^{n-2})
\]
may be chosen to fix the coattaching face $D^2\times S^{n-3}$ pointwise.
Hence the collapse on the next borrowed handle may be chosen to fix $J_i$.
\end{lemma}

\begin{proof}
Apply Cohen's product-collapse lemma \cite[Lemma~8.3]{Coh69} to
$D^2\searrow\{0\}$ and the factor $D^{n-2}$.  It gives the finite PL
collapse
\[
 D^2\times D^{n-2}\searrow
 (D^2\times S^{n-3})\cup(\{0\}\times D^{n-2})
\]
relative to its target, and therefore pointwise relative to the coattaching
face.  For $k_i^2$ that face contains precisely the protected trace $J_i$.
\end{proof}

\begin{lemma}[Relative based handle-slide and cancellation]
\label{lem:gt-relative-based-handle-cancellation}
For $n\ge6$, the handle choices in each borrowed block may be made so that,
after the old borrowed handle is carved out, the actual $2/3$-handle pair
is geometrically complementary in the sense of the based handle calculus of
\cite{RS72}.  Moreover, the handle-slide and cancellation may be transported
back to the literal triad $(A_i';F_i,F_{i+1})$ while fixing the retained old
handle chart and $J_i$ pointwise.  These choices may be made recursively,
one block ahead.
\end{lemma}

\begin{proof}
Use the based handles
\[
 H=h_i^1,
 \qquad A=h_i^2,
 \qquad K=k_i^2,
 \qquad L=k_i^3.
\]
Let $e$ be the relative $1$-handle generator and let $a,k$ be the
$2$-handle generators.  In right-module notation over the pre-surgery group
ring, write
\[
 \partial a=e u,
 \qquad
 \partial k=e v.
\]
The original geometric cancellation of $K/H$ makes $v$ a trivial unit.
Set $x=a_i$ and $y=a_{i+1}$.  The attaching word of $A$ is
\[
 w_A=y^{-2}x^{-1}yx,
\]
while that of $K$ is $w_K=y$.  After the old $a_i$-surgery, $x=1$ and
$w_A=y^{-1}$.  Write bars for the images under the induced quotient of the
based group ring after $x$ is killed.  Thus both $\bar u$ and $\bar v$ are
trivial units, with signs and group elements determined by the chosen
orientations, lifts, and basepaths.

In the original cancellation order the residual generator and the
$3$-handle boundary have the form
\[
 a'=a-kv^{-1}u,
 \qquad
 \partial L=a'\alpha,
\]
where $\alpha$ is a trivial unit.  Cancel $A/H$ instead and use the based
slide
\[
 \bar k'=\bar k-\bar a\,\bar u^{-1}\bar v.
\]
Then $\partial\bar k'=0$, and substitution gives
\[
 \partial\bar L
   =-\bar k'\,\bar v^{-1}\bar u\,\bar\alpha,
\]
again a single trivial unit.  At the word level, choose the orientations,
lifts, and basepaths compatibly with this basis change.  The $K$-word is
$y$ and the $A$-word is $y^{-1}$ after $x$ is killed, so the based slide
replaces the $K$-word by $yy^{-1}=1$, up to the basepath conjugacies already
encoded by the trivial units.  Thus the attaching loop of $K'$ is
null-homotopic.  This is the non-simply-connected based form of the Adding
Lemma described in \cite[p.~89]{RS72}; its piping arc records the required
group element.

We check the hypotheses of \cite[Corollary~6.20]{RS72}.  Write
$W'=W\cup K'\cup L$, where $W$ is the pre-pair product cobordism and $M_0$
is its incoming boundary.  Put
\[
 M_1=\partial W-M_0,
 \qquad
 M_2=\partial(W\cup K')-M_0.
\]
Thus $M_1$ is the outgoing boundary before $K'$ and $M_2$ is the
intermediate boundary after $K'$ and before $L$.  The corollary requires the
natural identifications
\[
 \pi_1(M_1)\cong\pi_1(M_2)\cong\pi_1(W).
\]
The first holds because $W$ is a product.  The null-homotopic attaching
loop of $K'$ shows that $W\hookrightarrow W\cup K'$ induces an isomorphism
on fundamental groups.  Viewed from $M_2$, the same handle has dual index
$n-2\ge4$, so
\[
 \pi_1(M_2)\xrightarrow{\cong}\pi_1(W\cup K')\cong\pi_1(W),
\]
which gives the second identification.  The incidence coefficient is a
trivial unit, so the proof of Corollary~6.20 uses the
non-simply-connected Whitney lemma to make $K'/L$ geometrically
complementary, after which the Cancellation Lemma
\cite[Lemma~6.4]{RS72} applies.  Its index range $2\le r\le n-4$, with
$r=2$, is exactly $n\ge6$.  The final outgoing boundary of $W'$ also has
the same fundamental group, since $L$ has index $3$ and the dual indices of
the pair are $n-3,n-2\ge3$.

It remains to preserve the literal markings.  Retain the notation of the
carving argument in \cite[pp.~283--284]{GT03},
\[
 E_i''=A_i''\cup_{\Delta_i}B_i''.
\]
The handles $H,A$ lie on the auxiliary $B_i''$ side, while $K,L$ are the
actual handles on the $A_i''$ side.  Canceling $H/A$ reparametrizes
$\Delta_i$.  The handle $K'$ is the same abstract actual handle $k_i^2$,
with its characteristic coordinates transported through the attachment
quotient; its ambient image is not asserted to remain fixed.  The auxiliary
handle-slide inserts no new handle into the actual block.  The based Adding Lemma just used, together with the attaching-map
isotopy of \cite[Lemma~6.1]{RS72}, gives an \emph{exchanged} handle
triad
\[
 (A_{i,\mathrm{ex}}'';F_{i,\mathrm{ex}},F_{i+1,\mathrm{ex}})
\]
and a homeomorphism of triads
\[
 \Phi_i:
 (A_{i,\mathrm{ex}}'';F_{i,\mathrm{ex}},F_{i+1,\mathrm{ex}})
 \longrightarrow
 (A_i'';F_i,F_{i+1}).
\]
Here the subscript $\mathrm{ex}$ denotes the handle presentation
obtained after the above handle slide and change of cancellation order.
Choose the future slide and Whitney $2$-spines and a closed coattaching
patch that will be the entire seam $J_i$ before freezing the Figure~5 map on
the newly created handle.  Make the spines disjoint from the old $2$-core
in the $(n-1)$-dimensional boundary; this is possible because
$2+2<n-1$.  Their $3$-dimensional traces in the block are disjoint from that
core because $2+3<n$.  Choose the coattaching patch in the unused part of
the surviving face, off all of these supports.  The attaching-map isotopies
are then performed relative to that whole patch; relative PL isotopy
extension realizes them by ambient homeomorphisms supported in its
complement.  Only afterward choose the retained old chart and the thin,
intermediate, and full regular neighborhoods sufficiently small.

The induction invariant is that, before the map on a newly created handle is
frozen for the next block, the next exchange support and its literal
coattaching seam have already been selected, and both the retained chart and
the whole seam lie in the fixed complement.  Thus $\Phi_i$ fixes the
retained old chart and the protected seam $J_i$ in the surviving
coattaching face.  Pull every post-slide chart back by $\Phi_i^{-1}$ and
conjugate every post-slide map accordingly.  The exchanged outgoing face is
thereby carried to the literal $F_{i+1}$ and the incoming face to the
literal $F_i$.  The homeomorphism need not fix all of $F_{i+1}$, and the
statement does not claim compatibility with an arbitrary previously frozen
Figure~5 chart.  The boundary epimorphism $G_{i+1}\twoheadrightarrow G_i$ is
unchanged up to the permitted basepath conjugacies.  This is a finite
collection of choices at each stage, so the construction continues by
countable induction without an infinite ambient composition.
\end{proof}

\begin{lemma}
\label{lem:gt-relative-common-celllike}
Let $A,B$ be subpolyhedra of a finite polyhedron $X$, put $J=A\cap B$, and
let $j:A\hookrightarrow X$ be inclusion.  Suppose
\[
 c:X\longrightarrow C,
 \qquad
 d=c|_A:A\longrightarrow C
\]
are onto cell-like maps to a finite polyhedron, and a map $r_B:B\to A$
satisfies
\begin{equation}\label{eq:relative-common-celllike}
 r_B|_J=\id_J,
 \qquad
 d r_B=c|_B.
\end{equation}
For every $\varepsilon>0$ and every metric on $C$, there are a retraction
$R:X\to A$, with $R|_B=r_B$, and a strong deformation retraction
\[
 H:X\times I\longrightarrow X,
 \qquad H_0=\id_X,
 \qquad H_1=jR,
 \qquad H_t|_A=\id_A,
\]
such that the $cH$-tracks and the $dRH$-tracks have diameter less than
$\varepsilon$.
\end{lemma}

\begin{proof}
After a common subdivision, apply Lacher's relative lifting lemma
\cite[Lemma~2.3]{Lac69} to $d:A\to C$, with complex pair $(X,A\cup B)$,
reference map $c$, and prescribed lift $\id_A\cup r_B$.  The first equality
in \eqref{eq:relative-common-celllike} makes the prescription a map, and the
second makes it an exact lift.  We obtain a retraction $R:X\to A$,
prescribed on $A\cup B$, for which $dR$ is arbitrarily close to $c$.

Finite polyhedra are locally equiconnected, so join $c$ to $dR$ by an
arbitrarily small homotopy $\Phi$ relative to $A$.  Apply Lacher's lemma a
second time to $c:X\to C$ and the pair
\[
 \bigl(X\times I,(X\times\{0,1\})\cup(A\times I)\bigr),
\]
prescribing $x$ at time zero, $R(x)$ at time one, and the stationary lift
$(a,t)\mapsto a$ on $A\times I$.  These are exact lifts of the
corresponding restrictions of $\Phi$.  The resulting lift is $H$.  The
first approximation, the smallness of $\Phi$, and the second approximation
give both track estimates.  Cell-like maps of finite polyhedra are
$UV^\infty$, so the dimension hypotheses of the lifting lemma are
satisfied.
\end{proof}

\begin{proposition}
\label{prop:gt-marked-twoface}
Suppose that, for all sufficiently large $i$, there are finite
subpolyhedra $C_i\subset F_i$, onto cell-like maps
\[
 d_i:F_i\longrightarrow C_i,
 \qquad
 c_i:A_i'\longrightarrow C_i,
\]
sections $s_i:C_i\hookrightarrow F_i$ with $d_is_i=\id_{C_i}$, and maps
$\bar f_i:C_{i+1}\to C_i$ such that
\begin{equation}\label{eq:gt-twoface-data}
 c_i|_{F_i}=d_i,
 \qquad
 c_i|_{F_{i+1}}=\bar f_i d_{i+1},
 \qquad
 d_i|_{J_i}=d_{i+1}|_{J_i}=\id_{J_i}.
\end{equation}
Then $M_*^n$ admits a finite-dimensional compact ANR
$\Z$-compactification.
\end{proposition}

\begin{proof}
For $x\in J_i$, the first and second equations in
\eqref{eq:gt-twoface-data} give $\bar f_i(x)=x$.  Put
\[
 \ell_i=s_i\bar f_i:C_{i+1}\longrightarrow F_i,
 \qquad
 f_i=\ell_i d_{i+1}:F_{i+1}\longrightarrow F_i.
\]
Then
\[
 f_i|_{J_i}=\id_{J_i},
 \qquad
 d_i f_i=c_i|_{F_{i+1}}.
\]
Apply \cref{lem:gt-relative-common-celllike} with
\[
 X=A_i',\qquad A=F_i,\qquad B=F_{i+1},\qquad r_B=f_i.
\]
It gives a retraction $R_i:A_i'\to F_i$, prescribed on both faces, and a
strong deformation retraction fixed on $F_i$.

The preceding bond restricts on $F_i$ as $f_{i-1}=\ell_{i-1}d_i$.
Uniform continuity of $\ell_{i-1}$ and of each of the finitely many earlier
stage maps now gives the required retained-coordinate estimate directly.
Namely, for $x\in A_i'$ and $j\le i-1$,
\[
 \rho_{j,i+1}H_i(x,t)
 =\rho_{j,i-1}\ell_{i-1}d_iR_iH_i(x,t),
\]
and the same identity at $t=0$ has $H_i(x,0)=x$.  After the earlier bonds
have been fixed, choose the relative Lacher tolerance simultaneously for the
finite family $j\le i-1$.  The $d_iR_iH_i$-track estimate then yields
\[
 \operatorname{dist}_j\!
 \left(\rho_{j,i+1}H_i(x,t),\rho_{j,i+1}(x)\right)<\frac1i.
\]
For $x\in K_i$ the homotopy is stationary, so the same estimate is automatic
there.  Paste the block maps to the identity on $K_i$ and proceed
inductively.  Now apply \cref{prop:gt-retained-control} with $a_i=i-1$ and
$\varepsilon_i=1/i$.
\end{proof}

\begin{proof}[Proof of \cref{thm:gt-marked-realization}]
Apply \cref{lem:gt-relative-based-handle-cancellation} and put the
complementary $K'/L$ pair in the coordinates of
\cref{lem:gt-complementary-pair-collapse} with
\[
 r=2,
 \qquad X=I^2,
 \qquad Z=I^{n-3},
 \qquad D=[1,3]\times Z.
\]
The protected seam lies in the surviving coattaching face of the pair, so
every local map below fixes $J_i$.  All charts and maps are conjugated by
the triad homeomorphism $\Phi_i$ of
\cref{lem:gt-relative-based-handle-cancellation}; consequently the notation
below refers to the literal actual faces.

We first choose the incoming Figure~5 map.  In a characteristic handle
$X\times D$, choose a concentric disk $D_0\subset\Int D$.  Cohen's
product-collapse lemma \cite[Lemma~8.3]{Coh69} gives a collapsible
retraction
\[
 \epsilon:X\times D\longrightarrow
 (X\times D-\Int D_0)\cup(\{x_0\}\times D_0).
\]
Apply \cref{lem:gt-one-sided-prism} radially across the outer annulus and
extend over $\{x_0\}\times D_0$ to obtain
\[
 \delta:(X\times D-\Int D_0)\cup(\{x_0\}\times D_0)
 \longrightarrow
 (X\times\partial D)\cup(\{x_0\}\times D).
\]
Thus $d=\delta\epsilon$ is the coattaching-relative Figure~5 retraction.
Both factors are collapsible and fix its target, so
\cite[Lemma~8.6]{Coh69} makes $d$ a collapsible retraction of $K$ onto $P$.
Use this prescribed $d$ in
\cref{lem:gt-complementary-pair-collapse}; henceforth $\varrho_i$ and $e_i$
denote the resulting maps, with the literal identity
\[
 \varrho_i|_K=e_i d.
\]
Transport this ordered pair of maps through every marked characteristic
chart.  On the handle $k_{i-1}^2$ denote the transported copies by
$\epsilon_i$ and $\delta_i$, and set $d_i=\delta_i\epsilon_i$; on the newly
created $k_i^2$ the same composite is $d_{i+1}$.  The thin, intermediate,
and full copies are chosen in this one chart when the handle is created.
This is a choice made one block ahead, not a factorization asserted for an
arbitrary previously fixed collapse.

Let
\[
 q_i:A_i'\longrightarrow
 Y_i=A_i''\cup\widetilde b_{i-1}^{\,n-2}
\]
be the first carving map of \cite[pp.~283--284]{GT03}, where
$\widetilde b_{i-1}^{\,n-2}$ is the retained thickened copy of the old belt
disk.  Choose it as the finite PL collapse retraction supplied by the
thin-handle shelling, fixing $Y_i$ pointwise.  Its incoming restriction is
$\epsilon_i$, while the retained full/intermediate shell carries
$\delta_i$; extend both by the identity on the complementary part of
$\Gamma_i$.  More explicitly, triangulate the thin, intermediate, and full
handle copies simultaneously.  Use \cite[Lemma~8.3]{Coh69} on the inner
copy, extend the specified elementary collapses through the radial shells by
\cref{lem:gt-one-sided-prism}, and collapse the remaining product collar
relative to its attaching face.  The old-carving support is disjoint from
the new complementary-pair ball, so every step extends over that ball by
the identity.  Cohen's composition lemma \cite[Lemma~8.6]{Coh69} gives the
stated collapsible PL retraction $q_i$ with the prescribed restrictions.

We now construct the map on the entire actual block and verify both face
restrictions.  Write the product part of the exchanged block as
$W_i^0=M_i\times I$, where $M_i$ is its incoming face, and put
$M_{i,1}=M_i\times\{1\}$.  Under these characteristic coordinates, write
\[
   \mathbb B_i=X\times[-1,3]\times Z,\qquad
   \mathbb K_i=X\times D,
\]
and let
\[
   \mathbb P_i
   =(X\times\partial D)\cup(\{x_0\}\times D)
   \subset \mathbb K_i
\]
be the Figure~5 target.  The attaching and final faces of
$\mathbb B_i$ are
\[
   Q_i=(\partial X\times[-1,3]\times Z)
       \cup(X\times\{-1\}\times Z),
\]
and
\[
   \mathbb R_i=
   X\times\bigl(([-1,3]\times\partial Z)\cup(\{3\}\times Z)\bigr).
\]  
Attach $\mathbb B_i$ to $M_{i,1}$ along $Q_i$.  Put
\[
 S_i^0=M_{i,1}-\Int_{M_{i,1}}Q_i
\]
and let $S_i$ be $S_i^0$ together with the retained old belt disk, chosen
disjoint from $\mathbb B_i$.  Let $\widehat W_i$ be $W_i^0$ together with
the retained old outer-shell and belt piece.  After the first carving,
\[
 Y_i=\widehat W_i\cup_{Q_i}\mathbb B_i,
 \qquad
 C_i=S_i\cup\mathbb R_i.
\]
Under the transported marking, the last space is the literal core
$\Gamma_i'\cup b_{i-1}^{n-2}\subset F_i$, and
\begin{equation}\label{eq:gt-pair-overlap}
 S_i\cap\mathbb B_i
 =\partial Q_i
 =Q_i\cap\mathbb R_i
 =\partial\mathbb R_i
 \subset\mathbb R_i.
\end{equation}
Use the same simultaneous triangulation.  Apply
\cref{lem:gt-one-sided-prism} to the specified elementary-collapse sequence
for $\delta_i$, thereby extending it over the retained old shell, and
perform the standard product-collar collapse to $M_{i,1}$ relative to
$Q_i$.  These maps agree on the common shell; their finite composition is
collapsible by \cite[Lemma~8.6]{Coh69}.  This gives a collapsible retraction
\[
 w_i:\widehat W_i\longrightarrow S_i\cup Q_i
\]
relative to $Q_i$.  Define
\[
 \widehat w_i=w_i\cup\id_{\mathbb B_i},
 \qquad
 \widehat\varrho_i=\id_{S_i}\cup\varrho_i,
\]
where $\varrho_i:\mathbb B_i\to\mathbb R_i$ is the pair retraction from
\cref{lem:gt-complementary-pair-collapse}.  The equalities in
\eqref{eq:gt-pair-overlap} make these extensions agree on their overlaps, so
the pasted maps are literal.  Each fiber of either extension is a singleton
or a fiber of the original collapsible retraction, so both extensions are
collapsible.  The same is true of $q_i$; hence Cohen's composition lemma
\cite[Lemma~8.6]{Coh69} shows that
\[
 c_i=\widehat\varrho_i\widehat w_i q_i:A_i'\longrightarrow C_i
\]
is a collapsible retraction.  On the incoming face the three restrictions
are $\epsilon_i$, $\delta_i$, and the identity, respectively.  Therefore
\[
 c_i|_{F_i}=\delta_i\epsilon_i=d_i
\]
as a literal equality.  Each factor fixes $C_i=S_i\cup\mathbb R_i$; hence
$c_i$ is literally a retraction onto $C_i$.

For the outgoing restriction, recall $D=[1,3]\times Z$ and decompose
\[
 F_{i+1}=\mathbb K_i\cup G_i,
 \qquad
 G_i=\Gamma_{i+1}-\Int_{\Gamma_{i+1}}(\partial X\times D).
\]
The overlaps are
\begin{equation}\label{eq:gt-outgoing-overlap}
 G_i\cap\mathbb P_i
 =G_i\cap\mathbb K_i
 =\partial X\times\partial D.
\end{equation}
The full next map $d_{i+1}$ is the local retraction
$d:\mathbb K_i\to\mathbb P_i$ glued to the identity on $G_i$; this gluing
is relative to $J_i$ by \cref{lem:gt-coattaching-figure5}.  In the composite
defining $c_i$, all stages before this local $d$ are the identity on the new
pair ball.  Let $h_i$ denote the single global composite after the local
$d$.  Then
\[
 h_i|_{\mathbb P_i}=e_i,
 \qquad
 g_i=c_i|_{G_i}.
\]
Because $d$ fixes $\mathbb P_i$, these maps agree on
\eqref{eq:gt-outgoing-overlap}; hence they glue to a PL map
\[
 \bar f_i:C_{i+1}=G_i\cup\mathbb P_i\longrightarrow C_i.
\]
On $\mathbb K_i$ one has
\[
 c_i=h_i d=e_i d=\bar f_i d_{i+1},
\]
whereas on $G_i$ one has $d_{i+1}=\id$ and
$\bar f_i=g_i=c_i|_{G_i}$.  Thus
\[
 c_i|_{F_{i+1}}=\bar f_i d_{i+1}
\]
literally.  Every local and pasted stage fixes the protected coattaching
face.  More explicitly, the protected-support choice and the prospective
definition of $d_i$ give $d_i|_{J_i}=\id_{J_i}$, while
\cref{lem:gt-coattaching-figure5} gives
$d_{i+1}|_{J_i}=\id_{J_i}$.

Freeze the transported chart on the actual $k_i^2$ and repeat.  The
countable induction consists only of finite choices at each stage.  The maps
just constructed satisfy \eqref{eq:gt-twoface-data}, so
\cref{prop:gt-marked-twoface} completes the proof.
\end{proof}

\begin{proof}[Proof of \cref{thm:boundaryless-example}]
Choose the marked Guilbault--Tinsley realization supplied by
\cref{thm:gt-marked-realization}.  The underlying actual blocks, their
frontiers, and their boundary epimorphisms have not changed.  Hence
\cite[Theorems~1.3 and~4.3]{GT03} gives a connected one-ended open PL
$n$-manifold whose end pro-fundamental group is semistable but not perfectly
semistable.  It is therefore not pseudo-collarable.  The preceding theorem
gives its finite-dimensional compact ANR $\Z$-compactification.

Finally, the $\Z$-boundary of any $\Z$-compactification of this manifold is
non-ANR.  Indeed, an ANR $\Z$-boundary would force stability by
\cref{prop:anr-expansion-restrictions}, whereas the Guilbault--Tinsley end
system is not even perfectly semistable.
\end{proof}

\section{Problems related to aspherical manifolds}
\label{sec:aspherical}

\subsection{Weak Borel and stability conjectures}

\begin{conjecture}[Weak Borel conjecture]\label{conj:weak-borel}
Closed aspherical manifolds with isomorphic fundamental groups have
homeomorphic universal covers.
\end{conjecture}

\begin{conjecture}[Aspherical manifold stability conjecture]
\label{conj:stability}
The universal cover of every closed aspherical manifold is inward tame.
\end{conjecture}

Both conjectures hold in dimensions at most $3$.  The weak Borel conjecture
was formulated in \cite{AG99}.  The stability conjecture is stronger than
Mihalik's semistability conjecture \cite{Mi83}, because inward tameness
implies semistability at infinity.  This comparison does not, however,
produce an implication from the stability conjecture to the weak Borel
conjecture.  Such an implication requires controlled compactifications, as
explained next.

\subsection{Controlled compactifications and boundary swapping}

\begin{definition}\label{def:controlled-compactification}
Let $(X,d)$ be a proper metric space and let
$(\overline X,\overline d)$ be a metrizable compactification.  It is a
\emph{controlled compactification} if, for every $\epsilon>0$ and $R>0$,
there is a compact set $C\subseteq X$ such that
\[
\operatorname{diam}_{\overline d}(B_d(x,R))<\epsilon
\]
whenever $B_d(x,R)\cap C=\varnothing$.
\end{definition}

This is the metric form of the nullity condition used in
\cite{GM19}.  The two-point compactification
$\mathbb R\sqcup\{-\infty,+\infty\}$ is controlled for the Euclidean metric.
By contrast, the square compactification of Euclidean $\mathbb R^2$ is not
controlled for the Euclidean metric.  These are statements about the chosen
metric and the chosen compactification; the mere existence of an unrelated
CAT(0) metric does not establish controlledness for the specified
metric--compactification pair.

A group action on a metric space is \emph{geometric} if it is proper,
cocompact, and isometric.  Following Bestvina \cite{Bes96}, we use the
following definition.

\begin{definition}\label{def:Z-structure}
A $\Z$-structure on a group $G$ is a pair $(\overline X,Z)$ such that
\begin{enumerate}[(1)]
\item $\overline X$ is a compact AR;
\item $Z$ is a $\Z$-set in $\overline X$;
\item $X=\overline X-Z$ is a proper metric space on which $G$ acts
geometrically; and
\item for every compact $K\subseteq X$ and every open cover $\mathcal U$ of
$\overline X$, all but finitely many $G$-translates of $K$ lie in members of
$\mathcal U$.
\end{enumerate}
The set $Z$ is then a $\Z$-boundary for $G$.
\end{definition}

For a geometric action, controlledness in
\cref{def:controlled-compactification} implies the nullity condition; see
\cite[Proposition~6.3 and Lemma~6.4]{GM19}.

\begin{theorem}[Boundary Swapping Theorem \cite{GM19}]
\label{thm:boundary-swapping}
Suppose that $G$ acts geometrically on proper metric ARs $X$ and $Y$.  If
$Y$ has a $\Z$-structure $(\overline Y,Z)$ for $G$, then the same compactum
$Z$ can be attached to $X$ to give a $\Z$-structure for $G$.
\end{theorem}

\begin{proposition}\label{prop:controlled-implies-weak-borel}
Let $M^n$ and $N^n$, $n\ge5$, be closed aspherical manifolds with a chosen
isomorphism $\pi_1(M)\cong\pi_1(N)=G$.  If $\widetilde M$ admits a controlled
$\Z$-compactification with respect to a proper $G$-invariant metric for
which the deck action is geometric, then $\widetilde M$ and $\widetilde N$
are homeomorphic.
\end{proposition}

\begin{proof}
The contractible manifolds $\widetilde M$ and $\widetilde N$ are proper
metric ARs, and $G$ acts geometrically on each by deck transformations for
the lifted metrics.  By \cite[Proposition~6.3 and Lemma~6.4]{GM19}, the
controlled $\Z$-compactification of $\widetilde M$ gives a $\Z$-structure
for $G$.  Boundary swapping attaches the same $\Z$-boundary to
$\widetilde N$.  The two contractible open $n$-manifolds therefore have
homeomorphic $\Z$-boundaries, so \cite[Theorem~12]{AG99} gives
$\widetilde M\cong\widetilde N$.
\end{proof}

\begin{corollary}
If every universal cover of a closed aspherical $n$-manifold, $n\ge5$,
admits a controlled $\Z$-compactification for its lifted metric, then the
weak Borel conjecture holds in dimension $n$.
\end{corollary}

An ordinary $\Z$-compactification supplied by
\cref{thm:retraction-intro} need not be controlled.  Thus vertical
topological control and metric nullity must not be conflated.  In particular,
\cref{prop:controlled-implies-weak-borel} does not presently follow from the
aspherical manifold stability conjecture.

\begin{question}\label{q:controlled-universal-cover}
Does the universal cover of every closed aspherical manifold admit a
controlled $\Z$-compactification for a lifted proper metric?
\end{question}

\subsection{Homotopy actions and controlled domination}

For completeness, we record the part of the transfer-reducibility framework
that follows from a finite-dimensional $\Z$-compactification.  Let $S$ be a
finite subset of a group $G$ containing the identity.

\begin{definition}[Homotopy $S$-action]
A homotopy $S$-action on a space $X$ consists of maps
$\varphi_g:X\to X$ for $g\in S$ and homotopies
\[
H_{g,h}:\varphi_g\circ\varphi_h\simeq\varphi_{gh}
\]
whenever $g,h,gh\in S$, with $\varphi_e=\id_X$ and the evident constant
identity homotopy.  The associated sets
$S^{n}_{\varphi,H}(g,x)$ and the notion of an $S$-long cover are those of
\cite[Definition~1.4]{BL12}.
\end{definition}

\begin{definition}\label{def:F-cover}
A \emph{family} $\mathcal F$ of subgroups of $G$ is closed under conjugation
and passage to subgroups.  If $Y$ is a $G$-space, a subset $U\subseteq Y$ is
an $\mathcal F$-subset if
\begin{enumerate}[(1)]
\item for every $g\in G$, either $gU=U$ or $gU\cap U=\varnothing$; and
\item $G_U=\{g\in G:gU=U\}$ belongs to $\mathcal F$.
\end{enumerate}
An open $\mathcal F$-cover is a $G$-invariant collection $\mathcal U$ of
open $\mathcal F$-subsets with $Y=\bigcup_{U\in\mathcal U}U$.
\end{definition}

\begin{definition}[Controlled $N$-domination]
A compact metric space $X$ is \emph{controlled $N$-dominated} if, for every
$\epsilon>0$, there are a finite CW complex $K$ of dimension at most $N$,
maps $i:X\to K$ and $p:K\to X$, and a homotopy from $p\circ i$ to $\id_X$
whose point tracks have diameter at most $\epsilon$.
\end{definition}

\begin{proposition}\label{prop:controlled-domination}
Let $M$ be a contractible manifold with a finite-dimensional compact metric
$\Z$-compactification $\overline M=M\sqcup Z$.  Then $\overline M$ is a
contractible compact ANR and is controlled $N$-dominated for some finite
$N$.  If a group $G$ acts on $M$, then every finite $S\subseteq G$ containing
the identity induces a homotopy $S$-action on $\overline M$.
\end{proposition}

\begin{proof}
The inclusion $M\hookrightarrow\overline M$ is a homotopy equivalence because
$Z$ is a $\Z$-set.  Hence $\overline M$ is contractible.  A
$\Z$-compactification of an ANR is an ANR.

Since $\overline M$ is finite-dimensional, it embeds as a closed subset of
some $\mathbb R^N$ and is a retract of an open neighborhood $U$.  Choose a
finite polyhedral neighborhood $K$ with
$\overline M\subseteq K\subseteq U$.  If
$i:\overline M\hookrightarrow K$ is inclusion and $p:K\to\overline M$ is
the neighborhood retraction, then $p\circ i=\id_{\overline M}$.  Thus the
controlled domination has constant tracks and works for every $\epsilon$.

Let $j:M\hookrightarrow\overline M$.  The $\Z$-set homotopy supplies a
homotopy inverse $r:\overline M\to M$ for $j$.  For $g\in S-\{e\}$, set
$\varphi_g=j\circ g\circ r$ and set $\varphi_e=\id$.  Inserting the
homotopies $rj\simeq\id_M$ and $jr\simeq\id_{\overline M}$ gives the required
$H_{g,h}$ whenever $g,h,gh\in S$.
\end{proof}

Transfer reducibility requires, in addition, uniformly bounded-dimensional
open $\mathcal F$-covers of $G\times\overline M$ that are $S$-long
\cite[Definition~1.8]{BL12}.  Here the genuine $G$-action on the product is
left translation on the first factor,
\[
k\cdot(g,x)=(kg,x),
\]
while $\overline M$ carries only the homotopy $S$-action constructed above.
Neither a $\Z$-compactification nor
\cref{prop:controlled-domination} automatically provides those covers.

\begin{question}
For which families $\mathcal F$ does the fundamental group of a closed
aspherical manifold with finite-dimensional $\Z$-compactifiable universal
cover admit the $S$-long covers required for transfer reducibility?
\end{question}

\end{document}